\documentclass[reqno,10pt]{amsart}
\usepackage{amsthm}
\usepackage{amsmath}
\usepackage{amssymb}
\usepackage{amsfonts}
\usepackage{mathrsfs}
\usepackage{todonotes}
\usepackage{verbatim}
\usepackage[bookmarksnumbered,colorlinks,urlcolor=blue]{hyperref}
\usepackage{graphicx}
\usepackage{graphics}
\usepackage{float}
\usepackage{enumerate}
\usepackage{ulem}
\def\bb#1\eb{\textcolor{blue}
{#1}} %
\def\br#1\er{\textcolor{red}
{#1}} %
\newcommand{\qcd}{\begin{flushright} $\Box$ \end{flushright}}

   \def\br#1\er{\textcolor{red}{#1}} %
      \def\bb#1\eb{\textcolor{blue}{#1}} %

\title[]{Path-lifting properties of the exponential map with applications to geodesics}

\author[I.P. Costa e Silva]{Ivan P. Costa e Silva}
\address{Department of Mathematics, Universidade Federal de Santa Catarina,
\hfill\break\indent 88.040-900
Florianópolis-SC, Brazil.}
\email{pontual.ivan@ufsc.br}

\author[J.L. Flores]{Jos\'e L. Flores}
\address{Departamento de \'Algebra, Geometr\'{\i}a y Topolog\'{\i}a,  Universidad de M\'alaga
\hfill\break\indent
Facultad de Ciencias, Campus Universitario de Teatinos,
\hfill\break\indent 29071 M\'alaga, Spain}
\email{floresj@uma.es}

\author[K.P.R. Honorato]{Kledilson P. R. Honorato}
\address{Department of Mathematics, Universidade Federal de Santa Catarina,
	\hfill\break\indent 88.040-900
	Florianópolis-SC, Brazil.}
\email{kledilson.honorato@posgrad.ufsc.br}

\begin{document}
\newtheorem{thm}{Theorem}[section]
\newtheorem{prop}[thm]{Proposition}
\newtheorem{lemma}[thm]{Lemma}
\newtheorem{cor}[thm]{Corollary}
\theoremstyle{definition}
\newtheorem{defi}[thm]{Definition}
\newtheorem{notation}[thm]{Notation}
\newtheorem{exe}[thm]{Example}
\newtheorem{conj}[thm]{Conjecture}
\newtheorem{prob}[thm]{Problem}
\newtheorem{rem}[thm]{Remark}
\newtheorem{conv}[thm]{Convention}
\newtheorem{crit}[thm]{Criterion}
\newtheorem{claim}[thm]{Claim}

\newcommand{\ben}{\begin{enumerate}}
\newcommand{\een}{\end{enumerate}}

\newcommand{\bit}{\begin{itemize}}
\newcommand{\eit}{\end{itemize}}

\begin{abstract}
We revisit certain path-lifting and path-continuation properties of abstract maps as described in the work of F. Browder and R. Rheindboldt in 1950-1960s, and apply their elegant theory to exponential maps. We obtain thereby a number of novel results of existence and multiplicity of geodesics joining any two points of a connected affine manifold, as well as {\it causal} geodesics connecting any two causally related points on a Lorentzian manifold. These results include a generalization of the well-known Hadamard-Cartan theorem of Riemannian geometry to the affine manifold context, as well as a new version of the so-called Lorentzian Hadamard-Cartan theorem using weaker assumptions than global hyperbolicity and timelike 1-connectedness required in the extant version. We also include a general discription of {\it pseudoconvexity and disprisonment} of broad classes of geodesics in terms of suitable restrictions of the exponential map. The latter description sheds further light on the the relation between pseudoconvexity and disprisonment of a given such class on the one hand, and geodesic connectedness by members of that class on the other. 

\end{abstract}

\maketitle
\section{Introduction}\label{sec0}

Countless mathematical questions boil down to the abstract problem of studying when a suitably defined map $F:X\rightarrow Y$ is surjective; and that issue is in turn connected with the existence of solutions $x \in X$ of (possibly nonlinear) equations of the form 
  \[
  F(x)=y.
  \]
 This viewpoint goes back in systematic form to J. Hadamard \cite{hadamard} in the early 1910s, and a number of elegant topological tools have been developed to address aspects of that general problem. Of special interest to us here is the analysis of {\it path-lifting} and {\it path-continuation} properties of certain maps by F. Browder and W. Rheinboldt \cite{bowder,rein} in the 1950s and 1960s. The resulting framework is so flexible that it has been adapted to a vast array of different problems (see, e.g., \cite{gutu,jaramillo} and references therein for a recent review with a number of applications). 
 
 Our goal in this paper is to apply aspects of the Rheinboldt-Browder theory to {\it exponential maps} to revisit, in a broad, unified perspective, both the problem of {\it geodesic connectedness} of affine manifolds and the problem of {\it causal geodesic connectedness} for Lorentzian manifolds. 

Let $(M,\nabla)$ be an {\it affine manifold}, i.e., $M$ is a smooth $m$-dimensional manifold ($m\geq 2$) and $\nabla$ is an affine conection thereon. The problem of {\it ($\nabla$-)geodesic connectedness} is as follows: if $M$ is connected, can any two points be connected via a geodesic segment? This problem is especially important for semi-Riemannian manifolds $(M,g)$ with their Levi-Civita connection.

Of course, trivial examples - say, like Euclidean spaces minus a point - suffice to show that the answer is negative in general. Thus, the question is actually under which geometrically natural hypotheses does geodesic connectedness occur. This problem has a completely satisfactory, well-known general solution for a (connected) {\it geodesically complete} Riemannian manifold via the Hopf-Rinow theorem, in which the existence of {\it minimal} geodesics connecting any two points of $M$ is established. Moreover, such completeness is automatic for {\it compact} Riemannian manifolds, and hence geodesic connectedness always occur for compact Riemannian manifolds. In addition, a classic result by Morse \cite[Thm. 13.3, p. 239]{M} (with an important addition by Serre \cite{serre}) establishes that any two {\it non-conjugate} points on a complete, non-contractible Riemannian manifold can be connected by infinitely many geometrically distinct geodesics. Even in the absence of geodesic completeness, however, the geodesic connectedness
of Riemannian manifolds is fairly well-understood \cite{bartolo,miguel}. 

Each and every one of the pleasant features pointed out above changes dramatically for more general affine manifolds, and even for indefinite semi-Riemannian manifolds, where no analogue of the Hopf-Rinow theorem exists. A famous example by Bates \cite{bates} has shown that even complete {\it and} compact affine manifolds may fail to be geodesically connected. Even if one only considers the more restricted (but greatly important) class of {\it Lorentzian manifolds}, de Sitter and anti-de Sitter spacetimes provide well-known examples of geodesically complete, maximally symmetric Lorentz manifolds which are not geodesically connected. The underlying manifolds in these examples, however, are not compact. Yet, compact analytic Lorentz tori
can still be found which also fail to be geodesically connected \cite[Example  7]{beemparker}. These examples, in turn, are not geodesically complete. Indeed, to the best of our knowledge it remains an open problem to ascertain whether a connected Lorentzian manifold which is both compact and geodesically complete is geodesically connected. (Recall that even for Lorentzian manifolds compactness does not imply geodesic completeness.)

Many techniques have been devised to study the problem of geodesic connectedness. These include powerful variational tools \cite{benci,annamaria,masiello}, existence of convex functions \cite{AlexKarr}, methods from group theory \cite{calabi}, and topological techniques based
on Brouwer's degree (see \cite{tesejose} and references therein), and more direct but restricted methods based on a (partial) integration of the geodesic equations (see, for instance, \cite{godel}). 

In the special case of Lorentzian manifolds, a related problem is that of {\it causal geodesic connectedness}: if two points in a Lorentzian manifold $(M,g)$ can be connected by a causal curve, then when can they be connected by a {\it causal geodesic}? This problem is of particular interest in gravitational physics, since the geometry of Lorentzian manifolds is at the heart of the general-relativistic description of gravity. In this context, signals between events in spacetime travel along causal curves therein, and causal geodesics then represent either light rays or free-falling observers under no influence other than that of the gravitational field. 

Just like in the Riemannian case, a general and completely satisfactory solution to the problem of causal geodesic connectedness occurs only when the Lorentzian manifold $(M,g)$ is a {\it globally hyperbolic spacetime} 
 \cite{avez,beem,seifert}. Indeed, the solution here is a perfect Lorentzian analogue of the complete Riemannian case: if $p,q \in M$ are connected by a causal curve, then they are connected by a {\it maximal} causal geodesic. The idea of the proof is purely variational in nature: global hyperbolicity implies compactness of a suitably defined space of causal curves, and a concomitant upper semicontinuity of the Lorentzian length functional does the rest.

Our approach here has been partially inspired by the work of J. Beem and P. Parker \cite{beem0,bp,beemparker}, who give conditions for geodesic connectedness in terms of {\it disprisonment} and {\it pseudoconvexity} of families of geodesics. While variational methods usually need functionals built out of a metric and hence apply mostly to the semi-Riemannian setting, these geometric conditions also make sense in the purely affine case. In addition, while geodesic completeness is often needed in other approaches, in the Beem-Parker approach (as well as in our own results here) the completeness assumption is largely irrelevant. This presents a great advantage especially in the Lorentzian case, where most physically interesting examples are not complete. 

 The key observation we make is that the Beem-Parker approach can be taken a step further by reinterpreting the purely geometric conditions of pseudoconvexity and disprisonment - even when they are valid only for restricted classes of geodesics such as timelike or null geodesics on an indefinite semi-Riemannian manifold - in terms of {\it analytic} properties of the exponential map $\exp :\mathcal{D}\subset TM \rightarrow M$ on the underlying affine manifold. This idea was hinted at in \cite{beemparker}, and systematically used and generalized in a recent paper \cite{AMPA} by two of us (IPCS and JLF) to get a geodesic connectedness result. 
 We greatly expand on that approach here, and in particular show how pseudoconvexity and disprisonment for certain subfamilies of geodesics - such as the causal geodesics in the semi-Riemannian context - can be reinterpreted as the condition making (suitable restrictions of) the exponential map {\it proper} (i.e., inverse images of compact sets are compact). As it turns out, properness of maps is primary among a number of properties ensuring the required path-lifting/path-continuation properties in the Rheinboldt-Browder theory, a fact which can then be used to our advantage. (The approaches in \cite{beemparker,AMPA}, however, dealt indiscriminately with {\it all} geodesics.) 
 
  Apart from properness of the exponential maps, we also discuss here a few weaker continuation conditions that still yield the existence of connecting geodesic/causal geodesic where appropriate, and in some circumstances also provide {\it multiplicity} statements, i.e., information about how many (inequivalent) geodesic segments of a certain type are there connecting two points. In fact, we are able to provide here both an affine generalization and a Lorentzian version of the well-know {\it Hadamard-Cartan theorem} from Riemannian geometry (cf. Cor. \ref{corclave1} and Thm. \ref{lorentzianhadamard} below).

 The rest of the paper is organized as follows. 
 
 In section \ref{sec1} we discuss some of the key results and notation we shall use in the subsequent sections. We present a short discussion of the Rheinboldt-Browder theory for the benefit of those unfamiliar with the theme. Our presentation is focused only on those results germane to us here, referring to the original references for proofs.
 
 In section \ref{actualsec3} we introduce an abstract description of pseudoconvexity and disprisonment properties for rather broad classes $\mathcal{C}$ of geodesics, exploring its announced relation with properness of the exponential map. This is done via a suitably defined set of vectors - which we call $C$ there - and the takeaway is summarized in Thm. \ref{prop1}. 
 
 In section \ref{sec2} we apply the results and notions of section \ref{sec1} to show that a suitable restriction of the exponential map is a covering map. That, in turn, is used to derive a novel geodesic connectedness result (Thm. \ref{thmclave1}) for affine manifolds which expands on the main theorem in \cite{AMPA}. As indicated above, a specialized application of the latter result yields a true generalization of the Hadamard-Cartan theorem (Cor. \ref{corclave1}), as well as a few subsidiary results. 
 
 In section \ref{sec3} we narrow our focus to Lorentzian manifolds, introducing a number of conditions based on the Rheinboldt-Browder theory of section \ref{sec1}, culminating in Thm. \ref{thmclave3} and its corollaries.
 
 Finally, in section \ref{sec4} we present a new Lorentzian version of the Hadamard-Cartan theorem (cf. Thm. \ref{lorentzianhadamard}) which assumes neither global hyperbolicity nor 1-connectedness like the previous version in Ref. \cite{beem}, and discuss briefly the placement of causally pseudoconvex and disprisoning spacetimes in the causal ladder \cite{minguzzi,ladder} (cf. Cor. \ref{causalladderplace}).

\section{Technical preliminaries}\label{sec1}

Our purpose in this section is twofold. First, to lay down most of the terminology and notation we will use for the rest of the paper. 

Second, we wish to describe a slightly abstract description of path-lifting  and continuation properties of maps between Hausdorff topological spaces, as put forth especially in Refs. \cite{bowder,jaramillo,rein}. This approach provides some generalizations of standard facts in the theory of covering spaces. Its importance for us lies in that it gives simple but elegant sufficient conditions for a map to be onto, and in some situations indeed a covering map. These notions will later be applied to the exponential map. 

Some features of this general framework have been adapted to, and found many applications in, analytic and geometric problems (see, e.g., \cite{jaramillo} , and also \cite{gutu} for a recent review with abundant references). Since its details may be unfamiliar to some geometers, we've decided to include a short account of those key points relevant to this paper. 


Throughout this section, let $X,Y$ denote arbitrarily fixed Hausdorff topological spaces, and $F:X\rightarrow Y$ is a continuous map between them. Following \cite{rein}, 
$\mathcal{P}(X)$ [resp. $\mathcal{P}(Y)$] denotes the set of all continuous paths $\alpha:[0,1]\rightarrow X[Y]$. (However, when applying the underlying topological theory specifically to smooth manifolds later on, we shall denote the latter by $M,N$ rather than $X,Y$ to emphasize the context.)

\subsection{Path-lifting and extension properties of maps}\label{subsecpathlift}

Some of the definitions and proofs below can be found in Refs. \cite{bowder,AMPA,jaramillo,rein}. Accordingly, we shall refer the reader to the relevant portions of those papers whenever appropriate, and present proofs only when they are not (to the best of our knowledge) to be found elsewhere. 

\begin{defi} \cite[Def. 2.3]{rein} \label{pathliftdefi}
Let a set $P\subset \mathcal{P}(Y)$ be given. We say that $F:X\rightarrow Y$ has the {\it path-lifting property for $P$} if for any $\alpha \in P$ and any $x_0 \in F^{-1}(\alpha(0))$ there exists a path $\overline{\alpha} \in \mathcal{P}(X)$ such that 
$$\overline{\alpha}(0) = x_0 \hbox{  and   } F\circ \overline{\alpha } = \alpha .$$
Following standard terminology, we refer to $\overline{\alpha}$ as a {\it lift of $\alpha$ through $F$} ({\it starting at $x_0$}). 
\end{defi}

Of course, if $F:X\rightarrow Y$ has the path-lifting property for a set $P\subset \mathcal{P}(Y)$, then it also has the path-lifting property for any $P'\subset P$. 

It is a basic property in the theory of covering spaces that if $F:X\rightarrow Y$ is a covering map, then it has the path-lifting property for $P = \mathcal{P}(Y)$, and hence to any of its subsets\footnote{Recall that $F$ is a {\it covering map} if it is onto and any $y\in Y$ has an {\it evenly covered} neighborhood $V \ni y$, i.e., $F^{-1}(V)\subset X$ is a disjoint union of open sets restricted to each of which $F$ is a homeomorphism onto $V$.}. The distinction of Rheinboldt's approach is both that this property can be weakened to proper subsets of $\mathcal{P}(Y)$, and also that he gave a very important characterization of path-lifting in terms of the following key definition. 

\begin{defi} \cite[Def. 2.2]{rein} \label{contpropdefi}
Let a set $P\subset \mathcal{P}(Y)$ be given. We say that $F:X\rightarrow Y$ has the {\it continuation property for $P$} if for any $\alpha \in P$, and any continuous function $\overline{\alpha}:[0,b)\subset [0,1] \rightarrow X$ such that 
$$F\circ \overline{\alpha } = \alpha |_{[0,b)},$$
there exists a sequence $(t_k)_{k\in \mathbb{N}}$ in $[0,b)$ converging to $b$ for which the sequence $(\overline{\alpha}(t_k))_{k\in \mathbb{N}}$ converges in $X$. 
\end{defi}

In these terms, the following important characterization of the path-lifting property, valid for {\it local homeomorphisms}, can be given. 

\begin{prop}\cite[Thm. 2.4]{rein} \label{reinthm}
Let $P\subset \mathcal{P}(Y)$. A local homeomorphism $F:X\rightarrow Y$ has the path-lifting property for $P$ if and only if it has the continuation property for $P$. 
\end{prop}
\qcd
It is in particular clear that a {\it covering} map $F:X\rightarrow Y $ has {\it both} the path-lifting property {\it and} the path-continuation property for {\it any} $P\subset \mathcal{P}(Y)$. 

The importance of the above discussion becomes apparent if we consider a general (maybe nonlinear) equation of the general form 
\begin{equation}\label{genericeq}F(x) =y_0\end{equation}
for a given $y_0 \in Y$. In this case the {\it existence} of solutions $x\in X$ of (\ref{genericeq}) is related to the {\it surjectivity} of $F$, and it becomes important to give suitable criteria for this. Following \cite{jaramillo} we say that $Y$ is {\it $P$-connected} for $P\subset \mathcal{P}(Y)$ if any two points on $Y$ can be connected by a path in $P$\footnote{Actually, the authors of \cite{jaramillo} add the requirement that for any path $\alpha:[0,1]\rightarrow Y$ in $P$, the reverse path $\alpha^{-1}:t\in [0,1]\mapsto \alpha(1-t) \in Y$ also belongs to $P$. In that paper, this is a technical condition we shall not use here, and therefore we omit it.}. It is clear that we have the following nice criterion for the surjectivity of $F$: {\it if $Y$ is $P$-connected and $F:X\rightarrow Y$ has the path-lifting property for $P$, then $F$ is onto.} By virtue of Prop. \ref{reinthm}, in the case when $F:X\rightarrow Y$ is a local homeomorphism the latter criterion can be restated as follows: {\it if $Y$ is $P$-connected and the local homeomorphism $F:X\rightarrow Y$ has the continuation property for $P$, then $F$ is onto}.

\subsection{Continuation properties for manifolds}\label{subsecweak}

When dealing with smooth manifolds\footnote{Here and hereafter, {\it (smooth) manifold} means a real, $C^{\infty}$, finite-dimensional, Hausdorff, second-countable manifold.} - our almost exclusive interest later in this paper - stronger results can be obtained. Let $N,M$ be smooth manifolds of the same (finite) dimension and let $F:N\rightarrow M$ be a smooth map. It is convenient to consider the collection of curves given by 
\begin{equation}\label{piece}P =\mathcal{P}^{\infty}(M):= \{\alpha \in \mathcal{P}(M) \, : \, \alpha \hbox{ is piecewise smooth}\}.\end{equation}

We shall make use of the following well-known result.
\begin{lemma}\cite[Thm. 7.28]{oneill}\label{oneilllemma1}
Let $N,M$ be semi-Riemannian manifolds, and let $F:N\rightarrow M$ be a local isometry. Assume that $M$ is connected and that, given any geodesic $\sigma:[0,1] \rightarrow M$ and any point $p\in N$ such that $F(p)=\sigma(0)$, there exists a lift of $\sigma$ through $F$ starting at $p$. Then $F$ is a smooth  (semi-Riemannian) covering map.
\end{lemma}
\qcd
We then have\footnote{Alternatively, it is also possible at the cost of few more definitions to establish Proposition \ref{contisenough} by adapting \cite[Thm. 2.6]{jaramillo}, but we shall not pursue the details here.}
\begin{prop}\label{contisenough}
Let $N,M$ be smooth manifolds with $M$ connected, and let $F:N\rightarrow M$ be a local diffeomorphism with the continuation property for $\mathcal{P}^{\infty}(M)$ [as defined in (\ref{piece})]. Then $F$ is a smooth covering map.
\end{prop}
\noindent {\it Proof.} Let $h$ be any Riemann metric on $M$, and consider the pullback $\overline{h}:=F^{\ast}h$. Since $F$ is a local diffeomorphism, this is also a well-defined Riemannian metric on $N$ with respect to which $F:(N,\overline{h})\rightarrow (M,h)$ becomes a local isometry. Any $h$-geodesic $\sigma:[0,1] \rightarrow M$ is an element of $\mathcal{P}^{\infty}(M)$ and for any point $p\in N$ such that $F(p)=\sigma(0)$, the continuation property of $F$ for $\mathcal{P}^{\infty}(M)$ ensures that there exists a lift $\overline{\sigma}:[0,1] \rightarrow N$ of $\sigma$ through $F$ starting at $p$ by Prop. \ref{reinthm}; since $F$ is a local isometry this lift is uniquely defined and an $\overline{h}$-geodesic. Lemma \ref{oneilllemma1} now yields that $F$ is a smooth covering map, as desired.
\qcd

Recall that a map $F:X\rightarrow Y$ between topological spaces $X,Y$ is {\it proper} if the inverse images of compact sets in $Y$ are compact in $X$. The next result is an immediate application of \cite[Thm. 7]{bowder} by recalling that a smooth proper map between manifolds is in particular {\it closed}, i.e., it maps closed sets onto closed sets (cf. also \cite[Prop. 4.6]{lee}). 
\begin{prop}\label{resultinbowder}
Let $N,M$ be smooth manifolds of same (finite) dimension and let $F:N\rightarrow M$ be a smooth local diffeomorphism. If $M$ is connected and $F$ is proper, then $F$ is a covering map, and moreover, the inverse image under $F$ of each point of $M$ is a {\it finite} set. 
\end{prop}
\qcd
 
\section{Pseudoconvexity and disprisonment via exponential maps}\label{actualsec3}
Throughout the rest of this paper we take $(M,\nabla)$ to be a connected affine manifold with dimension $dim\, M\geq 2$ and $\nabla$ is an affine connection thereon\footnote{In dealing with geodesics, there is no loss of generality in assuming $\nabla$ is torsion-free. The domain $I\subset \mathbb{R}$ of a geodesic $\gamma:I \rightarrow M$ may be an open, closed or half-closed non-empty interval unless specified, but in case $I=[a,b]$ is compact, we often refer to $\gamma$ as a geodesic {\it segment}, with {\it endpoints} $\gamma(a),\gamma(b)$.}. 

As indicated in the Introduction, we wish to consider the problem of geodesic connectedness by geodesics in a suitable, maybe proper subcollection $\mathcal{C}$ of all the geodesics on $(M,\nabla)$. Such a restriction to a suitable subclass is especially important when we particularize to a semi-Riemannian manifold $(M,g)$ (with its Levi-Civita connection) of indefinite signature. We assume the reader is familiar with at least the basic aspects of these. Our main reference for facts about semi-Riemannian manifolds is O'Neill's classic text \cite{oneill}.


Recall that on a semi-Riemannian manifold $(M,g)$ we have three distinguished classes of vectors $v\in TM$: {\it timelike} when $g(v,v)<0$, {\it null} or {\it lightlike} if $v\neq 0$ and $g(v,v)=0$, and {\it spacelike} if either $v=0$ or $g(v,v)>0$. The vector $v$ is said to be {\it causal} (or {\it nonspacelike}) if it is either null or timelike. This general classification is usually referred to as the {\it causal character} of vectors. This notion can be extended to smooth curves and vector fields on $M$ in a natural way: the curve/vector field is said to be timelike [resp. null, causal, spacelike] if the tangent vector at each point where it is defined is of the corresponding character. General curves may not have a definite causal character but geodesics always do, and we refer to them as timelike/null/causal/spacelike according to which causal character they have.

For each $p\in M$ we denote by $\mathcal{T}_p \subset T_pM$ the set of timelike vectors at $p$. In the special case when $(M,g)$ is a {\it Lorentzian} manifold $\mathcal{T}_p$ is the disjoint union of two connected open convex cones called {\it timecones}, and a continuous choice of one of these components - called the {\it future timecone} - throughout $(M,g)$ is called a {\it time orientation} on the Lorentzian manifold $(M,g)$. If such a choice is fixed the Lorentzian manifold is said to be {\it time-oriented}. A connected, time-oriented Lorentzian manifold is called a {\it spacetime}. Again, we assume the reader is familiar with basic aspects of these. The reader may consult the standard textbook \cite{beem} as needed. 

We denote by $\mathcal{C}_p$ the closure of $\mathcal{T}_p$ in $T_pM$. Note that $0_p \in \mathcal{C}_p$, and that $\mathcal{C}_p \setminus \{0_p\}$ coincides with the set of causal vectors in $T_pM$. 

It is convenient to set up first a slightly abstract framework to describe a suitable fixed (nonempty) collection $\mathcal{C}$ of geodesics on $(M,\nabla)$ in terms of the exponential map thereon. The general description can then be particularized for the cases of interest later. 

We denote by $\exp :\mathcal{D}\subset TM \rightarrow M$ the exponential map of $(M,\nabla)$, with $\mathcal{D}$ being its maximal domain. (In particular, $\mathcal{D}=TM$ if and only if $(M,\nabla)$ is geodesically complete.) In addition, $\mathcal{Z}\subset TM$ shall denote the range of the zero section of the tangent bundle $TM$ and $\pi:TM\rightarrow M$ its canonical projection; recall that $\mathcal{D}$ is open in $TM$, contains  $\mathcal{Z}$, and for any $v\in \mathcal{D}$ and $t\in [0,1]$ we have $t\cdot v\in \mathcal{D}$. Finally, if $p\in M$, $\mathcal{D}_p:=\mathcal{D}\cap T_pM$ and $\exp _p:= \exp |_{\mathcal{D}_p}$. 

Recall that a collection $\mathcal{C}$ of (non-constant) geodesics on $(M,\nabla)$ is 
\begin{itemize}
\item[a)] {\it pseudoconvex} if for any compact set $K\subset M$ there exists a compact set $K^{\ast}\subset M$ such that any segment of a geodesic in $\mathcal{C}$ with endpoints in $K$ is entirely contained in $K^{\ast}$;
\item[b)] {\it disprisoning} if for any given maximal extension $\gamma:(a,b)\rightarrow M$ of a geodesic in $\mathcal{C}$ ($-\infty\leq a<b\leq +\infty$), and any $t_0\in (a,b)$, neither $\overline{\gamma[t_0,b)}$ nor $\overline{\gamma (a,t_0]}$ is compact. If $\mathcal{C}$ is not disprisoning, then it is said to be {\it imprisoning}. 
\end{itemize}

If we take $\mathcal{C}$ to be the collection of {\it all} non-constant geodesics in $(M,\nabla)$ and that class is pseudoconvex and disprisoning, then $(M,\nabla)$ itself is said to be a ({\it geodesically}) {\it pseudoconvex} and {\it disprisoning} affine manifold.

For a semi-Riemannian manifold $(M,g)$, however, it is also interesting to consider separately the cases when $\mathcal{C}$ is the collection of all {\it null} (or {\it lightlike}) [resp. {\it timelike}, {\it causal}] geodesics. If, say, the collection of all causal geodesics is pseudoconvex/disprisoning, then $(M,g)$ itself is said to be {\it causally pseudoconvex}/{\it causally disprisoning}, respectively.

For convenience, we fix thoughout a {\it complete} auxiliary Riemannian metric $h$ on $M$. We shall denote its distance by $d_h$ and the $h$-norm of $v\in TM$ by $|v|_h$.  Let $S\subset M$ be any set. Then $\mathbb{S}^1_hS$ will denote the subset of $\pi^{-1}(S)\subset TM$ of vectors with unit $h$-norm. It is easy to check that if $S$ is compact in $M$, then $\mathbb{S}^1_hS$ is compact in $TM$, a fact that we will use often. 

In order to give a technically convenient description of the collection $\mathcal{C}$, fix a non-empty set $C\subset \mathcal{D}$ with the following properties.
\begin{itemize}
    \item[P1)] For any $v\in C$, there exists $v' \in C\setminus \mathcal{Z}$ with $\pi(v') = \pi(v)$. 
    \item[P2)] For any $v\in C\setminus \mathcal{Z}$, $(\mathbb{R}\cdot v) \cap \mathcal{D} \subset C$. In particular, if $v\in C$, then $0\cdot v \in C$. 
\end{itemize}
Then, define
\[
C_1:= \{v/|v|_h\in \mathbb{S}^1_hM \, : \, v \in C\setminus \mathcal{Z}\},
\]
and for each $w\in C_1$, consider the interval
\[
J_w:=\{t\in \mathbb{R} \, :\, t\cdot w \in \mathcal{D}\}.
\]
(Note that $0\in J_w$.) Finally, define the class of geodesics
\begin{equation}\label{theclass}
\mathcal{C}:= \{\gamma : t\in I \mapsto \exp (t\cdot w) \in M\, : \, I\subset J_w \mbox{ non-degenerate interval and }w\in C_1 \}. 
\end{equation}
About these definitions, the following comments are in order.
\begin{itemize}
    \item[a)] For each $w\in C_1$, property $(P2)$ of $C$ implies that $J_w$ is an open interval, and $t\in J_w \mapsto \exp(t\cdot w) \in M$ is a maximally extended geodesic. 
    \item[b)] $\mathcal{C}$ is non-empty (because $C\neq \emptyset$), and because of $(P1)$ and the requirement that intervals are nondegenerate, there are no constant geodesics in $\mathcal{C}$.
    \item[c)] Given any geodesic in $\mathcal{C}$, any of its extensions are in $\mathcal{C}$. Furthermore, given any geodesic in $\mathcal{C}$, and its maximal extension $\gamma:(a,b)\rightarrow M$, for any $t_0 \in (a,b)$ we have that $\gamma|_{(a,t_0]}$ and $\gamma|_{[t_0,b)}$ also belong to $\mathcal{C}$. 
\end{itemize}
\bigskip

\noindent Let us consider some concrete examples of this general setting. Not all of these will be considered later in this paper, but they will still serve to illustrate the fairly broad scope of the previous definitions.   

\begin{exe}\label{ex1}
For $(M,\nabla)$ arbitrary, take $C\equiv \mathcal{D}$. We claim that $\mathcal{C}$ coincides with the collection of {\it all nonconstant} geodesics in $(M,\nabla)$ (up to affine reparametrizations) in this case. First, note that  $C_1\equiv  \mathbb{S}^1_hM$ in this case. Let $\alpha :I\rightarrow M$ be any nonconstant geodesic, and pick any $t_0\in I$. Let $v:=\alpha'(t_0)$. Since $\alpha$ is non-constant, $v\neq 0$, so we define 
\[
w:= \frac{v}{|v|_h} \in \mathbb{S}^1_hM. 
\]
Now, we have, for $t\in I$,
\[
(t-t_0)\cdot v \in \mathcal{D} \mbox{ and } \alpha(t)= \exp((t-t_0)\cdot v) = \exp((t-t_0)|v|_h\cdot w) \Rightarrow (t-t_0)|v|_h \in J_w. 
\]
Thus, if we define 
\[
I_w:= \{(t-t_0)|v|_h \, : \, t\in I\}
\]
we conclude that $I_w\subset J_w$, so $\gamma: s\in I_w \mapsto \exp(s\cdot w) \in M$ is in $\mathcal{C}$ and is an affine reparametrization of $\alpha$ as desired. 
\end{exe}
\begin{exe}\label{ex2}
Again for arbitrary $(M,\nabla)$, pick any $p\in M$ and take $C\equiv \mathcal{D}_p$. Then $C\cap \mathcal{Z}=\{0_p:=0_{T_pM}\}$. Reasoning exactly as in the previous example, given a nonconstant geodesic $\alpha:I\rightarrow M$ with $\alpha(t_0)=p$ for some $t_0\in I$, we can show that up to an affine reparametrization this is a geodesic in $\mathcal{C}$. We conclude that $\mathcal{C}$ comprises all nonconstant geodesics (which can be extended to pass) through $p$ up to an affine reparametrization.
\end{exe}
\begin{exe}\label{ex3}
Let $(M,g)$ be a semi-Riemannian manifold, and let $S\subset M$ be a smooth semi-Riemannian (embedded) submanifold of codimension $>0$. Denote by $NS\subset TM$ its normal bundle, and put
\[
C= NS\cap \mathcal{D}.
\]
(So that $\exp_{\perp}:= \exp|_{C}$ is the normal exponential of $S$.) It is not difficult to check that in this case $\mathcal{C}$ is the collection of all nonconstant geodesics in $(M,g)$ which can be extended as a normal geodesic to $S$, up to affine reparametrization. 
\end{exe}
\begin{exe}\label{exe4}
Let $(M,g)$ be a semi-Riemannian manifold. Recall that for each $p \in M$ we denote by $\mathcal{T}_p$ the set of all timelike vectors $v\in T_pM$. We then put
\[
C= (\mathcal{T}_p\cup \{0_p\})\cap \mathcal{D}.
\]
Then (again up to affine reparametrizations and/or extensions) $\mathcal{C}$ is the collection of all timelike geodesics passing through $p$. Again, if $ \Lambda _p$ denotes the set of {\it lightlike} (nonzero) vectors in $T_pM$, then we may define 
\[
C=(\Lambda_p\cup \{0_p\})\cap \mathcal{D}
\]
and $\mathcal{C}$ is the collection of all lightlike geodesics through $p$. Finally, let $\mathcal{C}_p$ be the closure of $\mathcal{T}_p$ in $T_pM$. Note that $0_p \in \mathcal{C}_p$, and that $\mathcal{C}_p \setminus \{0_p\}$ coincides with the set of causal vectors in $T_pM$. For 
\[
C= \mathcal{C}_p\cap \mathcal{D}
\]
then (again up to reparametrizations/extensions) $\mathcal{C}$ will be the class of all causal geodesics through $p$. 
\end{exe}
We can now state our first main result. It generalizes results in \cite{beemparker} (cf. also \cite[Prop. 2.1(i)]{AMPA}), and gives a convenient analytic characterization of pseudoconvexity-cum-disprisonment for a broad range of geodesic families.

\begin{thm}\label{prop1}
If $\exp |_{C}$ is a proper map, then the collection $\mathcal{C}$ defined in (\ref{theclass}) is pseudoconvex and disprisoning. If $C\cap \mathcal{Z}$ is compact and $C$ is closed in $\mathcal{D}$, then the converse holds. 
\end{thm}
\noindent {\it Proof.} Suppose $\exp |_{C}:C\rightarrow M$ is proper. Assume, by way of contradiction, that  $\mathcal{C}$ is imprisoning.  Consider, in that case, an inextendible geodesic $\gamma:(a,b)\rightarrow M$ (with $-\infty\leq a<b\leq +\infty$) in $\mathcal{C}$ and some $t_0\in (a,b)$ for which $\overline{\gamma[t_0,b)}$ is compact (the case when $\gamma |_{(a,t_0]}$ is compact being entirely analogous). If we write $\gamma(t)= \exp(t\cdot w)$ for every $t\in [t_0,b)$ and some $w\in C_1$, we have
\[
t\cdot w \in \tilde{K}:= (\exp |_{C})^{-1}(\overline{\gamma[t_0,b)}), \quad \forall t\in [t_0,b),
\]
with $\tilde{K}$ compact in $C$, and hence also in $TM$. Let $p:=\pi(w)$, so $w\in T_pM$. Clearly $t\cdot w\in T_pM\cap \tilde{K}$, for every $t\in [t_0,b)$, and since $T_pM\cap \tilde{K}$ is compact in $T_pM$, it is also $h_p$-bounded, so we necessarily have $b<+\infty$. But then we also have $b\cdot w\in \tilde{K}\subset C$, whence $b\in J_w$, i.e., $\gamma$ would be right-extendible, contradicting the assumed inextendibility of $\gamma$. This establishes that $\mathcal{C}$ is disprisoning.

Assume now that  $\mathcal{C}$ is not pseudoconvex. In particular, $M$ cannot be compact. Pick any $q\in M$, and for each $n\in \mathbb{N}$ let
\[
\overline{B}_n:= \{p \in M \, : \, d_h(p,q)\leq n\}. 
\]
Each $\overline{B}_n$ is compact by the Hopf-Rinow theorem for $(M,h)$, since it is closed and $d_h$-bounded and $h$ is complete. The failure of pseudoconvexity also implies that there exists a compact set $K\subset M$ and a sequence $(\gamma_n:[a_n,b_n]\rightarrow M)_{n\in \mathbb{N}}$ of elements of $\mathcal{C}$ such that
\[
\gamma _n(a_n),\gamma _n(b_n)\in K \mbox{ and } \exists t_n\in [a_n,b_n] \mbox{ with } \gamma _n(t_n)\notin \overline{B}_n, \quad \forall n \in \mathbb{N}.
\]
Write 
\[
\gamma_n(t)=\exp (t\cdot w_n), \quad w_n \in C_1, t\in [a_n,b_n], \forall n \in \mathbb{N}. 
\]
Then, 
\[
a_n\cdot w_n,b_n\cdot w_n \in \tilde{K}:=  (\exp |_{C})^{-1}(K), \quad \forall n \in \mathbb{N}, 
\]
and again $\tilde{K}$ is compact in $C$ and hence in $TM$. In addition, for each $n\in \mathbb{N}$
\[
p_n:= \pi(w_n)=\pi(a_n\cdot w_n)= \pi(b_n\cdot w_n) \in \pi(\tilde{K}).
\]
The latter set is compact in $M$, hence so is $\mathbb{S}^1_h\pi(\tilde{K})\ni w_n$. We conclude that we may pick a subsequence $(\gamma_{n_i})_{i\in \mathbb{N}}$, $w\in \mathbb{S}^1_hM$ and $a,b \in \mathbb{R}$ such that $ w_{n_i} \rightarrow w$ in $TM$, $a_{n_i}\rightarrow a$ and $b_{n_i}\rightarrow b$ in $\mathbb{R}$. In particular $w\neq 0$ and $a\leq b$. Furthermore, since $\tilde{K}$ is closed in $C$ we also have $a\cdot w,b\cdot w\in \tilde{K} \subset C \subset \mathcal{D}$, and the fact that $a_{n_i}\leq t_{n_i}\leq b_{n_i}$ forces (again up to passing to subsequences) a convergence $t_{n_i}\rightarrow t_0\in [a,b]$. Now, as $a\cdot w,b\cdot w\in \mathcal{D}_{\pi(w)}$ we also have $[a,b]\cdot w \subset \mathcal{D}_{\pi(w)}$ as the latter set is star-shaped. This in turn implies that $t_0 \cdot w \in \mathcal{D}$, and thus
\[
\gamma_{n_i}(t_{n_i}) = \exp (t_{n_i}\cdot w_{n_i})\rightarrow \exp(t_0\cdot w). 
\]
In particular $(\gamma _{n_i}(t_{n_i}))$ would be $d_h$-bounded, which is impossible by its construction. This contradiction establishes that $\mathcal{C}$ is indeed pseudoconvex. 

To prove the converse, assume that  $\mathcal{C}$ is pseudoconvex and disprisoning, with $C\cap \mathcal{Z}$ compact and $C$ closed in $\mathcal{D}$. Let $K\subset M$ be a compact set and pick any sequence $(v_n) \subset  (\exp |_{C})^{-1}(K)$. We wish to show that this has a subsequence converging to $(\exp |_{C})^{-1}(K)$. If all but a finite number of terms in this sequence are zero vectors we of course do have such a subsequence, so passing to a subsequence if necessary we may assume that every $v_n\neq 0$, and define
\[
w_n:= \frac{v_n}{|v_n|_h}, \mbox{ and } a_n:= |v_n|_h>0, \quad \forall n\in \mathbb{N}.
\]
By construction each $w_n\in C_1$, and $[0,a_n] \subset J_{w_n}$. Let $p_n:= \pi(w_n)\equiv \pi(v_n)$. Since $0_{p_n}\in C\cap \mathcal{Z}$, which we are assuming to be compact in $\mathcal{D}$ (and hence in $TM$), we have that 
\[
p_n \in \pi(C\cap \mathcal{Z}), \quad \forall n\in \mathbb{N},
\]
and the latter set is compact in $M$. But then $w_n \in \mathbb{S}^1_h\pi(C\cap \mathcal{Z})$, therefore, up to passing to a subsequence we may assume that there exists $w\in \mathbb{S}^1_hM$ for which $w_{n}\rightarrow w$. Now define, for each $n\in \mathbb{N}$, the geodesic
\[
\gamma _n: t\in [0,a_n] \mapsto \exp(t\cdot w_n) \in M,
\]
which is in $\mathcal{C}$ by construction. Consider the compact set $\mathcal{K}:=  \pi(C\cap \mathcal{Z})\cup K\subset M$. Since $\gamma_n(0)=p_n,\gamma_n(a_n)=\exp(v_n) \in \mathcal{K}$ for each $n\in \mathbb{N}$, by pseudoconvexity we can choose a compact set $\mathcal{K}^{\ast}\subset M$ such that $\gamma_n[0,a_n] \subset \mathcal{K}^{\ast}, \forall n\in \mathbb{N}$. Let $\gamma:[0,b)\rightarrow M$ be the right-inextendible geodesic with $\gamma'(0)=w$. Then $t\cdot w\in \mathcal{D}$ for any $t\in [0,b)$ and 
\[
\gamma(t) = \exp(t\cdot w), \quad \forall t\in [0,b).
\]
We now claim that $(a_n)$ is bounded. For suppose not. Passing to a subsequence if needed, we may assume in that case that $a_n\rightarrow +\infty$. Fix $0<t<b$. Eventually $a_n>t$, so on the one hand
\[
t\cdot w_n = \left(\frac{t}{a_n}\right)\cdot v_n \in C, \quad \forall n\in \mathbb{N},
\]
since $0< t/a_n<1$. On the other hand 
\[
t\cdot w_n \rightarrow  t\cdot w\in \overline{C}\cap \mathcal{D}\equiv \overline{C}^{\mathcal{D}}\equiv C
\]
since $C$ is closed in $\mathcal{D}$. We conclude that $t\cdot w \in C$, and since $t>0$ we have $w\in C_1$. In addition $\gamma_n(t)\rightarrow \gamma(t)$, so $\gamma(t)\in \mathcal{K}^{\ast}$. We conclude that $\gamma \in \mathcal{C}$ and $\gamma[0,b)\subset \mathcal{K}^{\ast}$ which contradicts the fact that $\mathcal{C}$ is disprisoning. 

Therefore $(a_n)$ is bounded as claimed. Again up to passing to a subsequence we assume that $a_n\rightarrow a$. Let $v:=a\cdot w$. Then $v_n\rightarrow v\in \overline{C}$. Thus all that remains to be shown is that $v\in \mathcal{D}$. But if not, $b\leq a$, so for any $0<t<b\leq a$, eventually $a_n>t$ and an argument just like the preceding one would establish that $\gamma \in \mathcal{C}$ and $\gamma[0,b)\subset \mathcal{K}^{\ast}$, again a contradiction. Thus $a<b$, so $v\in \mathcal{D}\cap \overline{C}= \overline{C}^{\mathcal{D}}\equiv C$, and $\exp(v) \in K$, so that $v\in  (\exp |_{C})^{-1}(K)$ as desired.
\qcd 

Using the notation in Example \ref{exe4}, with the choice
\[
C_p= \mathcal{C}_p\cap \mathcal{D}
\]
we can describe the class of all causal geodesics on $(M,g)$ emanating from $p$. Since such $C_p$ is closed in $\mathcal{D}$ and $C_p\cap \mathcal{Z}=\{0_p\}$ is compact, Theorem \ref{prop1} has the following immediate consequence.

In order to avoid repetition, henceforth whenever we use the phrase ``causally pseudoconvex and disprisoning'', it is to be understood that the latter part means {\it causally} disprisoning. 

\begin{cor}\label{causalcorshort}
Let $(M,g)$ be a Lorentz manifold, and $p \in M$. The set of causal geodesics in $(M,g)$ passing through $p$ (up to extensions) is causally pseudoconvex and disprisoning if and only if $\exp_p|_{C_p}$ is proper. In particular, if $(M,g)$ is itself causally pseudoconvex and disprisoning, then $\exp_q|_{C_q}$ is proper for {\it every} $q \in M$. 
\end{cor}
\qcd
\begin{rem}\label{GH1}
With an eye towards physical applications it is often convenient that geometric properties of spacetimes be derived from its position in the so-called {\it causal ladder} \cite{ladder,minguzzi}. The highest rung on that ladder is occupied by the {\it globally hyperbolic} spacetimes. It is well-known (cf., e.g., \cite[Prop. 7.36]{beem}) that if a spacetime $(M,g)$ is globally hyperbolic, then it is causally pseudoconvex and disprisoning. The converse, however, is false. This is illustrated by the strip $\{(t,x) \in \mathbb{R}^2 \, : \, 0<x<1\}$ (with the restricted metric) in the Minkowski plane $(\mathbb{R}^2,-dt^2 +dx^2)$. 
\end{rem}
\section{Geodesic connectedness on affine manifolds}\label{sec2}

In this section we employ the abstract results in section \ref{sec1} to study the geodesic connectedness of a general (connected) affine manifold $(M,\nabla)$. 

Let $p\in M$. The image by $\exp _p$ of the set 
$$\mathcal{S}_p:= \{v\in \mathcal{D}_p \, : \, (d\exp )_v \hbox{ is singular}\}$$
will be denoted by 
$$Conj (p) := \exp_p(\mathcal{S}_p) $$
and coincides with the set of all conjugate points to $p$ along some geodesic. 

Two of the authors [IPCS and JLF] have obtained in \cite{AMPA} a sufficient condition for a connected affine manifold to be geodesically connected. We say that $(M,\nabla)$ is a {\it weakly Wiedersehen} manifold, or $WW$-manifold for short, if for any $p \in M$ 
\begin{itemize}
    \item[a)] $Conj(p)$ is closed and
    \item[b)] $M\setminus Conj(p)$ is connected.
\end{itemize}
(In Thm. \ref{thmclave1} below, however, this condition makes sense - and only needs to be applied - pointwise. Therefore it is natural to say that $(M,\nabla)$ is $WW$ {\it at (a point) $p\in M$} if conditions $(a),(b)$ hold at $p$.) 

In \cite{AMPA}, the following concept was also introduced. 
\begin{defi}\label{def1.1}
Let $p\in M$. We say that $\exp _p$ is {\it weakly proper} if for any continuous\footnote{In \cite{AMPA}, the requirement is actually that $\overline{\alpha}$ is only (regular and) piecewise smooth; however the difference is really a technicality, because in the concrete proofs in which the concept is used the lift will always be contained in the region of $T_pM$ where $\exp_p$ is a local diffeomorphism; hence, $\overline{\alpha}$ can be chosen as a piecewise smooth lift of a piecewise smooth curve.} curve $\overline{\alpha}:[0,a) \rightarrow T_pM$ ($0< a\leq +\infty$) such that $\overline{\alpha}[0,a) \subset \mathcal{D}_p $ and $\exp_p \circ \overline{\alpha} :[0,a)\rightarrow M $ is right-extendible, there exists a compact set $K\subset \mathcal{D}_p$ containing $\overline{\alpha}[0,a)$.
\end{defi}
Using these concepts, a geodesic connectedness result was proved therein (cf.\cite[Thm 1.2]{AMPA}).  
\begin{thm}\label{mainthminAMPA}
Let $M=(M,\nabla)$ be a connected WW affine manifold. If $\exp_p:\mathcal{D}_p \subset T_pM \rightarrow M$ is weakly proper for each $p \in M$, then $M$ is geodesically connected. 
\end{thm}
\qcd
We now seek to establish an improved version of the latter theorem. In order to apply the Rheinboldt-Browder theory we first reinterpret weak properness in terms of a continuation property for the exponential. 
\begin{prop}\label{thenewprop1}
The map $\exp_p:\mathcal{D}_p \rightarrow M$ is weakly proper if and only if it has the continuation property for $\mathcal{P}^{\infty}(M)$.
\end{prop}
\noindent {\it Proof.} Let $\alpha \in \mathcal{P}^{\infty}(M)$ and consider any continuous map $\overline{\alpha}:[0,b)\subset [0,1] \rightarrow \mathcal{D}_p$ such that $\exp_p\circ \overline{\alpha} = \alpha |_{[0,b)}$. Assume that $\exp_p$ is weakly proper: thus, $\overline{\alpha}[0,b) \subset K$ for some compact set $K\subset \mathcal{D}_p$. 

 Let $(t_k)_{k\in \mathbb{N}}$ be {\it any} sequence in $[0,b)$ converging to $b$. Since $\overline{\alpha}(t_k)\in K$ for any $k$ we can, up to passing to a subsequence, assume that $(\overline{\alpha}(t_k))_{k\in \mathbb{N}}$ converges in $\mathcal{D}_p$.
 
For the converse, let $\overline{\alpha}:[0,b)\rightarrow \mathcal{D}_p\subset T_p M$, $0<b<1$, be a continuous curve with $\alpha:=\exp_p\circ\overline{\alpha}:[0,b)\rightarrow M$ right-extendible.
Assume that $\overline{\alpha}[0,b)$ is not contained in any compact set of $\mathcal{D}_p$, so in particular $\exp _p$ is not weakly proper. We shall construct another curve violating the continuation property for $\exp_p$.

Since $\exp_p$ has the continuation property for $\mathcal{P}^{\infty}(M)$, there exists a sequence $(t_k)_{k\in \mathbb{N}}$ in $[0,b)$ converging to $b$ for which the sequence $(\overline{\alpha}(t_k))_{k\in \mathbb{N}}$ converges in $\mathcal{D}_p$. 
Let $(\mathcal{K}_k)_{k\in \mathbb{N}}$ be an increasing (i.e. satisfying $\mathcal{K}_k\subset\mathring{\mathcal{K}}_{k+1}$) sequence of compact sets of $\mathcal{D}_p$ such that, given any compact set $\mathcal{K}\subset \mathcal{D}_p$, there exists $k_0$ such that $\mathcal{K}\subset \mathring{\mathcal{K}}_k$ for any $k\geq k_0$. 
Since $\overline{\alpha}[0,b)$ is not contained in any $\mathcal{K}_k$, and yet $(\overline{\alpha}(t_k))_{k\in \mathbb{N}}$ converges in $\mathcal{D}_p$, we can assume the existence of sequences $(t^l_k)_{k\in \mathbb{N}}$ with
\[
t_1<t^1_1<t^1_2<t_2<t^1_3<t^2_3<t^2_4<t^1_4<t_3<t^1_5<t^2_5<t^3_5<t^3_6<t^2_6<t^1_6<t_4<\cdots
\]
such that
\[
\overline{\alpha}(t_k)\in \mathring{\mathcal{K}}_1\quad\forall k\qquad \hbox{and}\qquad\overline{\alpha}(t^l_k)\in {\rm Fr}(\mathcal{K}_{l})\quad \forall k\geq 2l-1.
\]
In fact, we can suppose that the limit of $(\overline{\alpha}(t_k))_{k\in \mathbb{N}}$, and indeed, each of its points, belongs to $\mathring{\mathcal{K}}_1$. In particular, the point $\overline{\alpha}(t_1)$ is located on $\mathcal{K}_1$; then, $\overline{\alpha}$ escapes from $\mathcal{K}_1$ and goes back to $\mathring{\mathcal{K}}_1$ at $t_2$; so, it necessarily passes twice through the boundary ${\rm Bd}(\mathcal{K}_1)$ at some values $(t_1<)t^1_1<t^1_2(<t_2)$. This behavior for $\overline{\alpha}$ is repeated from the values $t_2$ to $t_3$, but now escaping from both $\mathcal{K}_1$ and $\mathcal{K}_2$, hence passing twice through the corresponding boundaries ${\rm Bd}(\mathcal{K}_1)$, ${\rm Fr}(\mathcal{K}_2)$ at some values $(t_2<)t^1_3<t^2_3<t^2_4<t^1_4(<t_3)$; and so on.

Let $h$ be some auxiliary complete Riemannian metric on $M$. Since ${\rm Bd}(\mathcal{K}_k)$ is compact for all $k$ and $\exp_p$ is smooth, there exist sequences $(m^1_n)_{n\geq 1}\supset (m^2_n)_{n\geq 2}\supset\cdots\supset (m^l_n)_{n\geq l}\supset\cdots$, with $m^l_{l+1}=m^{l+1}_{l+1}$ for all $l$, and piecewise smooth curves $\overline{\sigma}_{n,l}$ in $\mathcal{D}_p$ connecting $\overline{\alpha}(t^l_{m^l_n})$ with $\overline{\alpha}(t^l_{m^l_{n+1}})$, such that
\[
{\rm length}_h(\exp_p\circ\overline{\sigma}_{n,l})<1/2^l\qquad\forall n,l.
\]
In particular,
$\overline{\sigma}_{l,l}$ is a curve in $\mathcal{D}_p$ connecting $\overline{\alpha}(t^l_{m^l_l})=\overline{\alpha}(t^{l}_{m^{l-1}_l})$ with $\overline{\alpha}(t^l_{m^l_{l+1}})$ such that
\begin{equation}\label{ju}
	{\rm length}_h(\exp_p\circ\overline{\sigma}_{l,l})<1/2^l\qquad\forall l.
\end{equation}   
Consider the curve $\overline{\sigma}:[0,b)\rightarrow \mathcal{D}_p$ obtained by the concatenation of (appropriate reparametrizations of) the following segments:\footnote{The term $\overline{\alpha}^{\pm 1}$ means that
we must consider $\overline{\alpha}$ if the extreme values of the interval domain are well-ordered, or $\overline{\alpha}$ otherwise.}

\[
\overline{\alpha}\mid_{[0,t^1_{m^1_1}]},\quad \overline{\sigma}_{1,1},\quad \overline{\alpha}^{\pm 1}\mid_{[t^1_{m^1_2},t^2_{m^1_2}]},\quad \overline{\sigma}_{2,2},\quad \overline{\alpha}^{\pm 1}\mid_{[t^2_{m^2_3},t^3_{m^2_3}]}\ldots 
\]

Let $(s_k)_{k\in \mathbb{N}}$ be an {\it arbitrary} sequence in $[0,b)$ converging to $b$. Then, $(\overline{\sigma}(s_k))_{k\in \mathbb{N}}$ does not converge in $\mathcal{D}_p$, since it escapes from any compact subset of $\mathcal{D}_p$. However, $(\exp_p\circ\overline{\sigma}(s_k))_{k\in \mathbb{N}}$ is a Cauchy sequence,
and thus, $\overline{\sigma}$ violates the continuation property. Indeed, in order to check the Cauchy character of  $(\exp_p\circ\overline{\sigma}(s_k))_{k\in \mathbb{N}}$ observe that, for each $k$, 
\[
(a_k:=)\exp_p\circ\overline{\sigma}(s_k)\in \left\{\begin{array}{l} 
\alpha^{\pm 1}([\hat{t}_k,\tilde{t}_k])\quad \hbox{with $\hat{t}_k:=t^{l_k}_{m^{l_k}_{l_k+1}}$, $\;\;\tilde{t}_k:=t^{l_k+1}_{m^{l_k}_{l_k+1}}$} \\ \hbox{or}
\\
{\rm Im}(\exp_p\circ\overline{\sigma}_{l_k,l_k})
\end{array}\right.
\]

Moreover, since $l_k\rightarrow\infty$ as $k\rightarrow\infty$, $t^l_m>t_k$ whenever $l\geq k$ and $t_k\rightarrow b$ as $k\rightarrow\infty$, we deduce that $\hat{t}_k, \tilde{t}_k\rightarrow b$ as $k\rightarrow\infty$.
Recall also the right-extendibility of $\alpha\mid_{[0,b)}$, the condition (\ref{ju}), and the fact that each segment $\exp_p\circ\overline{\sigma}_{l_k,l_k}$ has its extreme points on $\alpha$. Therefore,  
given two arbitrary elements $a_{k_i}$, $i=1,2$, there are three possibilities: (i) if 
$a_{k_i}\in \alpha^{\pm 1}([\hat{t}_{k_i},\tilde{t}_{k_i}])$, $i=1,2$,
the distance between them is small if $k_1$, $k_2$ are big, as a consequence of the right-extendibility of $\alpha\mid_{[0,b)}$ and the fact that $\hat{t}_{k_i}, \tilde{t}_{k_i}\rightarrow b$ as $k_i\rightarrow\infty$; (ii) if, say, $a_{k_1}\in \alpha^{\pm 1}([\hat{t}_{k_1},\tilde{t}_{k_1}])$ and $a_{k_2}\in {\rm Im}(\exp_p\circ\overline{\sigma}_{l_{k_2},l_{k_2}})$, then
\[
d_h(a_{k_1},a_{k_2})\leq d_h(a_{k_1},p)+d_h(p,a_{k_2})\quad\hbox{for $p$ some extreme of $\exp_p\circ\overline{\sigma}_{l_{k_2},l_{k_2}}$.}
\]
Thus, one reaches the same conclusion as before, since $d_h(a_{k_1},p)$ is small by the right-extendibility of $\alpha\mid_{[0,b)}$ and $d_h(p,a_{k_2})$ is also small due to condition (\ref{ju}).
Finally, (iii) if $a_{k_i}\in {\rm Im}(\exp_p\circ\overline{\sigma}_{l_{k_i},l_{k_i}})$, $i=1,2$, we repeat the argument above by applying the triangle inequality to four points, where now the two ``middle'' points are extreme points of $\exp_p\circ\overline{\sigma}_{l_{k_i},l_{k_i}}$, $i=1,2$, resp.
 
\qcd

Using the notation and results in subsection \ref{subsecweak}, we finally have
\begin{thm}\label{thmclave1}
Let $(M, \nabla)$ be an affine manifold and $p \in M$. Assume that
\begin{itemize}
\item[i)] $\exp _p$ has the continuation property for $\mathcal{P}^{\infty}(M)$ (conf. Eq.(\ref{piece})), (or, equivalently, that it is weakly proper), and
\item[ii)] $(M,\nabla)$ is WW at $p$.
\end{itemize}
Then the map $\varphi_p:=\exp_p|_{\mathcal{N}_p}:\mathcal{N}_p \rightarrow M\setminus Conj(p)$ is a smooth covering map, where 
$$\mathcal{N}_p := \exp _p^{-1} ( M\setminus Conj(p)).$$
In particular, for any $q\in M\setminus Conj(p)$ there exists at least one geodesic segment from $p$ to $q$. 
\end{thm}
\noindent {\it Proof.} First of all, $(ii)$ implies that $\mathcal{N}_p$ and $M\setminus Conj(p)$ are open (and by Sard's Theorem, nonempty) in $\mathcal{D}_p$ and $M$, respectively.  Furthermore,
$$v\in \mathcal{N}_p \Rightarrow \exp _p(v) \notin Conj(p) \Rightarrow v\notin \mathcal{S}_p,$$
and hence $\varphi_p$ is a local diffeomorphism by the Inverse Function Theorem. 

The proof will be accomplished if we can apply Proposition \ref{contisenough} to $\varphi_p$. Thus, all that remains to be shown is that $\varphi_p$ has the continuation property for $\mathcal{P}^{\infty}(M\setminus Conj(p))$. 

To that end, fix $v_0 \in \mathcal{N}_p$ and let $\alpha:[0,1] \rightarrow M\setminus Conj(p)$ with $\alpha(0) =\varphi_p(v_0) =\exp_p(v_0)$ be a path in $\mathcal{P}^{\infty}(M\setminus Conj(p))$, which can naturally be viewed as a subset of $\mathcal{P}^{\infty}(M)$. Let $\overline{\alpha}:[0,l)\subset [0,1]\rightarrow \mathcal{N}_p$ be a continuous curve starting at $v_0$ for which 
$$\varphi_p\circ \overline{\alpha} \equiv \exp _p \circ \overline{\alpha} = \alpha |_{[0,l)}.$$

Because of the continuation property of $\exp_p$, there exists a sequence $(t_k) \subset [0,l)$ with $t_k \rightarrow l$ for which $ \overline{\alpha}(t_k)$ converges in $\mathcal{D}_p$ to some $\hat{v} \in \mathcal{D}_p$, say. But by continuity, $\exp_p(\hat{v}) = \alpha(l) \in \alpha[0,1] \subset M\setminus Conj(p)$, so $\hat{v} \in \mathcal{N}_p$. This establishes the continuation property as desired, so the proof is complete. 

\qcd
\begin{rem}\label{rem1}
The following comments about Theorem \ref{thmclave1} and Proposition \ref{thenewprop1} are in order.
\begin{itemize}
\item[1)] There are some concrete situations where weak properness is known to apply (conf. \cite[Props. 2.6, 2.7]{AMPA}). Suppose $(M,g)$ is either a complete Riemannian manifold or a geodesically complete Lorentz manifold possessing a {\it parallel timelike} vector field $V$. Then for any $p\in M$ the exponential map $\exp _p:T_pM \rightarrow M$ is weakly proper, although not necessarily proper (e.g., when $M$ is compact). 
    \item[2)] It may well happen that $p \in Conj(p)$, i.e., that $p$ is {\it self-conjugate}. In that case, $0_p \notin \mathcal{N}_p$. This occurs, for example, for any point on the round sphere $\mathbb{S}^m$ for any $m\geq 2$. Note that all the assumptions in Thm. \ref{thmclave1} hold here at any point, since $\mathbb{S}^m$ is a complete Riemannian manifold (cf. previous item) and for any $p \in \mathbb{S}^m$ we have $Conj(p) = \{p,-p\}$. The covering map $\varphi_p$ is trivial in this case for any $m\geq 2$: although there are infinitely connected components of $\mathcal{N}_p$, each is diffeomorphic to $\mathbb{S}^m \setminus \{p,-p\}$; this is expected for $m>2$ as $\mathbb{S}^m\setminus \{p,-p\}$ is simply connected, but holds also for $m=2$, even though $\pi_1(\mathbb{S}^2\setminus \{p,-p\})=\mathbb{Z}$. Thm. \ref{thmclave1} correctly predicts that for any $q\in \mathbb{S}^m \setminus \{p,-p\}$ there are infinitely many geodesic segments from $p$ to $q$, although in this case they are of course all parts of a single complete closed geodesic passing through $p$ and $q$. However, since $p$ is self-conjugate, the theorem does {\it not} predict the infinitely many closed geodesic at $p$ itself. 
    \item[3)] The previous discussion may be contrasted with the $2$-dimensional flat (Riemannian) cylinder $C$. Here there are no conjugate points, and $\mathcal{N}_p = T_pC\simeq \mathbb{R}^2$, so $\varphi_p$ is actually a universal covering of infinite multiplicity. For any $p,q \in C$, Thm \ref{thmclave1} correctly predicts that there are infinitely many geodesic segments connecting $p$ and $q$, even for $p=q$, although again all segments are part of a single closed geodesic at $p$ in the latter case. 
\end{itemize}
\end{rem}
The following generalized version of Hadamard-Cartan Theorem is a consequence of Thm. \ref{thmclave1} and a discussion in \cite{beemparker}. 
\begin{cor}[Generalized Hadamard-Cartan]\label{corclave1}
Let $(M^m, \nabla)$ be a connected affine manifold and let $p \in M$. Assume that $Conj(p) =\emptyset$ and that $\exp_p$ has the continuation property for $\mathcal{P}^{\infty}(M)$ (or, equivalently, that it is weakly proper). Then the following statements hold. 
\begin{itemize}
    \item[i)] Either $\pi_1(M)$ is trivial or (countably) infinite. 
    \item[ii)] If $M$ is simply connected, then it is actually diffeomorphic to $\mathbb{R}^m$. In this case there is a unique (up to affine reparametrizations) geodesic segment connecting any $q\in M$ and $p$. 
    \item[iii)] If $M$ is not simply connected, then for each $q\in M$ there exist countably infinitely many\footnote{Such segments may again be all parts of a single closed geodesic - cf. Remark \ref{rem1}(2).} geodesic segments from $p$ to $q$. In particular, there exist infinitely many geodesic loops at $p$. 
\end{itemize}
\end{cor} 
\noindent {\it Proof.} Using the notation in Thm. \ref{thmclave1}, we note that since $Conj(p) =\emptyset$ we have $\mathcal{N}_p \equiv \mathcal{D}_p$. Since the latter set is open and star-shaped in $T_pM$, it is diffeomorphic to $\mathbb{R}^m$. Now, since the conditions of Thm. \ref{thmclave1} apply ($(i)$ occurring because of Prop. \ref{thenewprop1}), it follows that $\exp_p: \mathcal{D}_p \rightarrow M$ itself is a (universal) covering map. Thus $(i)$ follows from \cite[Lemma 8]{beemparker}. 

As for $(ii)$, if $M$ is simply connected then the covering $\exp_p$ is trivial, and since it is connected, it is actually a diffeomorphism between $\mathcal{D}_p \simeq \mathbb{R}^m$ and $M$. 

Finally, $(iii)$ holds because the multiplicity of the cover coincides with the cardinality of $\pi_1(M)$, which is infinite in this case by $(i)$.
\qcd

\section{Causal geodesic connectedness on Lorentzian manifolds}\label{sec3}

In this section and in the next one we shall particularize the ideas in the previous sections to {\it Lorentzian} manifolds (endowed with the corresponding Levi-Civita connection). 

Fix therefore, for the rest of this section, a Lorentzian manifold $(M^m,g)$ with $m\ge2$. 

Recall (cf. Example \ref{exe4}) that for each $p\in M$ we denote by $\mathcal{T}_p \subset T_pM$ the set of timelike vectors at $p$. Remember that $\mathcal{T}_p$ is the disjoint union of two connected open convex cones called {\it timecones}. 

We denote by $\mathcal{C}_p$ the closure of $\mathcal{T}_p$ in $T_pM$. Note that $0_p \in \mathcal{C}_p$, and that $\mathcal{C}_p \setminus \{0_p\}$ coincides with the set of causal vectors in $T_pM$. Again, $\mathcal{C}_p\setminus \{0_p\}$ has two connected components called {\it causal cones}. A piecewise smooth curve $\sigma:[a,b]\rightarrow M$ is said to be {\it timelike} [resp. {\it causal}] if its tangent vector $\sigma '(t) \in \mathcal{T}_{\sigma(t)}$ [resp. $\in \mathcal{C}_{\sigma(t)}\setminus \{0_{\sigma(t)}\}$] for any $t\in [a,b]$ and both lateral tangent vectors at a break are on the same component of the timecone [resp. causal cone] thereat.

As in section \ref{actualsec3}, let
\begin{equation}\label{causalBigC}
C_p:= \mathcal{C}_p\cap \mathcal{D}.
\end{equation}
Following standard notation, we write
\begin{eqnarray}
I(p) &=&\{q\in M \, : \, \hbox{$\exists$ piecewise smooth timelike segment connecting $p$ and $q$}\}, \nonumber \\
J(p)&=& \{q\in M \, : \, \hbox{$\exists$ piecewise smooth causal segment connecting $p$ and $q$}\}\cup\{p\}.\nonumber\end{eqnarray}
It is well-known that $I(p)$ is always open.

Although Theorem \ref{thmclave1} has a nice conclusion including multiplicity of geodesics as well as existence, as a Lorentz-geometric result it has the drawback that it requires a continuation property for all smooth curves, say via weak properness, while in studying causal geodesic connectedness one would expect that control only along causal curves should suffice. We now proceed to show that this is indeed the case, although for technical reasons we will have to make a few slight adaptation of the continuation notions introduced by Rheinboldt and discussed in section \ref{sec1}. 

\begin{defi}[Causal continuation property]\label{def1}
Let $p\in M$. We say that $\exp _p$ has the {\it causal continuation property} (CCP) if for any (piecewise smooth) causal curve $\sigma:[0,1] \rightarrow M$ with $\sigma(0)=p$, and for any continuous curve $\overline{\sigma}:[0,a)\subset [0,1] \rightarrow C_p$ [for $C_p$ defined in (\ref{causalBigC})] such that $\overline{\sigma}(0)=0_p$ and 
$$\exp_p\circ \overline{\sigma} = \sigma|_{[0,a)}$$
there exists a sequence $(t_k)_{k \in \mathbb{N}}\subset [0,a)$ with $t_k\rightarrow a$ for which $\overline{\sigma}(t_k)$ converges to some $\overline{x}\in \mathcal{D}_p$ (and thus $\overline{x} \in C_p$). 
\end{defi}

Recall that if $(\hat{M},\hat{g})$ is a Lorentzian manifold, a map $\phi:\hat{M}\rightarrow M$ is a {\it Lorentzian covering map} if it is a smooth covering map for which $\hat{g} = \phi^{\ast}g$; in particular it is a local isometry. 

\begin{lemma}\label{lemma1.1}
Let $\phi:(\hat{M},\hat{g})\rightarrow (M,g)$ be a Lorentzian covering map. 
\begin{itemize}
    \item[i)] If $(M,g)$ is causally pseudoconvex and disprisoning, then so is $(\hat{M},\hat{g})$. 
\item[ii)] For each $\hat{p} \in \hat{M}$, $\exp ^{\hat{g}} _{\hat{p}}$ has the CCP if and only if $\exp ^g_{\phi(\hat{p})}$ does.
\end{itemize}
\end{lemma}
\noindent {\it Proof.} To simplify notation, we shall indicate by a hat any quantity pertaining to $(\hat{M},\hat{g})$, and those without are understood to belong to $(M,g)$. Fix $\hat{p} \in \hat{M}$ and $p=\phi(\hat{p}) \in M$. 
Using covering lifting properties of $\phi$ and the fact it is a local isometry, it is not difficult to check that 
\begin{equation}\label{eq1.1}
d\phi _{\hat{p}}(\hat{\mathcal{D}}_{\hat{p}}) = \mathcal{D}_p \hbox{ and } d\phi _{\hat{p}}(\hat{\mathcal{C}}_{\hat{p}}) = \mathcal{C}_p \Longrightarrow d\phi_{\hat{p}}(\hat{C}_{\hat{p}}) = C_p,
\end{equation}
as well as
\begin{equation}\label{eq1.2}
\phi\circ \hat{\exp }_{\hat{p}} = \exp _p \circ d\phi _{\hat{p}}. 
\end{equation}
$(i)$\\
Assume $(M,g)$ is causally pseudoconvex and disprisoning. Then $\exp _p|_{C_p}$ is proper by Cor. \ref{causalcorshort}. Let $\hat{K}\subset \hat{M}$ be a compact subset, and let $(\hat{v}_k)\subset (\hat{\exp}_{\hat{p}}|_{\hat{\mathcal{C}}_{\hat{p}}})^{-1}(\hat{K})$ be any sequence. The compactness of $\hat{K}$ implies that up to passing to a subsequence, we can assume that
$$\hat{\exp}_{\hat{p}}|_{\hat{\mathcal{C}}_{\hat{p}}}(\hat{v}_k) \rightarrow \hat{q} \Rightarrow \phi\circ \hat{\exp }_{\hat{p}}|_{\hat{\mathcal{C}}_{\hat{p}}}(\hat{v}_k)  \rightarrow \phi(\hat{q})=:q \stackrel{(\ref{eq1.2})}{\Rightarrow} \exp_p(d\phi _{\hat{p}}(\hat{v}_k) ) \rightarrow q. $$
Pick any compact neighborhood $V\ni q$ in $M$. Since $\exp _{p}|_{C_p}$ is proper, $W:= (\exp _{p}|_{C_p})^{-1}(V)$ is compact in $C_p$, and eventually $(d\phi _{\hat{p}}(\hat{v}_k))\subset W$. Again up to passing to a subsequence, then, 
$$d\phi _{\hat{p}}(\hat{v}_k) \rightarrow v_0$$
for some $v_0 \in C_p$. But $d\phi_{\hat{p}}$ is an isomorphism, so by (\ref{eq1.1})
$$\hat{v}_k\rightarrow (d\phi_{\hat{p}})^{-1}(v_0) =: \hat{v}_0 \in \hat{C}_{\hat{p}}.$$
We conclude that $\hat{\exp}_{\hat{p}}|_{\hat{C}_{\hat{p}}}(\hat{v}_0)\equiv \hat{q} \in \hat{K}$, i.e., $\hat{v}_0 \in (\hat{\exp}_{\hat{p}}|_{\hat{\mathcal{C}}_{\hat{p}}})^{-1}(\hat{K})$, whence the compactness of the latter set follows. \\
$(ii)$\\
We show only that if $\hat{\exp } _{\hat{p}}$ has the CCP, then so does $\exp _{p}$, since the proof of the converse is similar.


Fix a (piecewise smooth) causal curve $\alpha:[0,1]\rightarrow M$ with $\alpha(0)=p$ and a continuous $\overline{\alpha}:[0,a)\subset [0,1] \rightarrow C_p$  such that $\overline{\alpha}(0)=0_p$ and $\alpha|_{[0,a)}=\exp_p \circ \overline{\alpha}$. 

Let $\beta:[0,1] \rightarrow \hat{M}$ be the unique lift of $\alpha$ through $\phi$ starting at $\hat{p}$. (It is in particular a causal curve in $(\hat{M},\hat{g})$ since $\phi$ is a local isometry.)  Let $\overline{\beta} := (d\phi _{\hat{p}})^{-1}\circ \overline{\alpha}$. Then $\overline{\beta}(0)=0_{\hat{p}}$, and by (\ref{eq1.1})
\begin{equation}\label{eq1.3}
\overline{\beta}[0,a) \subset \hat{\mathcal{D}}_{\hat{p}} \cap \hat{\mathcal{C}}_{\hat{p}}.
\end{equation}

Now, 
\begin{eqnarray}
\phi\circ \hat{\exp}_{\hat{p}}\circ \overline{\beta} &\stackrel{(\ref{eq1.2})}{=}& \exp_p \circ d\phi _{\hat{p}}\circ \overline{\beta} \nonumber\\
&=& \exp _p \circ \overline{\alpha} \equiv \alpha |_{[0,a)}.\nonumber 
\end{eqnarray}
By the uniqueness of the lift of $\alpha|_{[0,a)}$ through $\phi$ starting at $\hat{p}$ we must have 
\begin{equation}\label{neweq1.4}
 \hat{\exp}_{\hat{p}}  \circ \overline{\beta} = \beta |_{[0,a)}. 
\end{equation}
Eq. (\ref{neweq1.4}) and the CCP for $\hat{\exp}_{\hat{p}}$ means there exists a sequence $(t_k)_{k \in \mathbb{N}}\subset [0,a)$ with $t_k\rightarrow a$ for which $\overline{\beta}(t_k)$ converges to some $\hat{x} \in \hat{C}_{\hat{p}}$. Thus $\overline{\alpha}(t_k) \rightarrow d\phi_{\hat{p}}(\hat{x})=: x \in C_p$ and we are done.
\qcd 
We now give a natural sufficient condition for exponential maps at points having the CCP.
\begin{prop}\label{prop2}
If $(M,g)$ admits a Lorentzian covering $\phi:(\hat{M},\hat{g})\rightarrow (M,g)$ such that $(\hat{M},\hat{g})$ is causally pseudoconvex and disprisoning, then for any $p\in M$ the exponential map $\exp _p$ has the CCP.
\end{prop}
\noindent{\it Proof.} By Lemma \ref{lemma1.1}(ii), there is no loss of generality in assuming that $(M,g)$ is itself causally pseudoconvex and causally disprisoning; hence we do so for the rest of the proof.

Fix any (piecewise smooth) causal curve $\sigma:[0,1] \rightarrow M$ with $\sigma(0)=p$, and any continuous curve $\overline{\sigma}:[0,a)\subset [0,1] \rightarrow C_p$ [for $C_p$ defined in (\ref{causalBigC})] such that $\overline{\sigma}(0)=0_p$ and 
$$\exp_p\circ \overline{\sigma} = \sigma|_{[0,a)}.$$

Using Corollary \ref{causalcorshort}, $\exp_p|_{C_p}$ is proper, and therefore 
$$K:= (\exp_p|_{C_p})^{-1}(\alpha[0,1])$$
is compact in $C_p$ and $\overline{\alpha}[0,a) \subset K$. 

Consider {\it any} sequence $(t_k)_{k \in \mathbb{N}}\subset [0,a)$ with $t_k\rightarrow a$. Since the image of $\overline{\alpha}$ is contained in a compact set in $C_p$, it follows that $\overline{\sigma}(t_k)$ converges to some $\overline{x}\in C_p$ up to passing to a subsequence. This establishes the CCP as desired.

\qcd

\begin{rem}\label{ccpweaker}
{\it The converse of Proposition \ref{prop2} is false}. To see this, consider again the $2$-dimensional Minkowski spacetime $(\mathbb{R}^2, -dt^2+dx^2)$. If we remove the entire fourth quadrant $Q = \{(t,x) \, : \, t\leq 0, x\geq 0\}$, the spacetime $(M:= \mathbb{R}^2\setminus Q, -dt^2 + dx^2)$ is simply connected, so any Lorentzian covering is trivial (thus a global isometry if connected). It is not causally pseudoconvex, but for each $p\in M$, it is easy to check that $\exp_p$ has the CCP.
\end{rem}

The following theorem aims at giving sufficient conditions to ensure the existence of a {\it timelike}  geodesic from $p \in M$ to $q\in I(p)$. Here, we denote by $Conj_c(p)$ the set of conjugate points to $p$ along {\it causal} geodesics starting at $p$\footnote{Note that $Conj_c(p)$ is contained in, but is {\it not} necessarily equal to, $Conj(p)\cap J(p)$.}. This is enough for causal connectedness, because if $q \in J(p) \setminus I(p)$, then it is well-known (cf. \cite[Prop. 10.46]{oneill}) that there exists a {\it null} geodesic segment connecting $p$ and $q$. 

\begin{thm}\label{thmclave3}
Let $(M, g)$ be a Lorentz manifold and $p \in M$; assume that $\exp _p$ has the CCP.
Let $q\in I(p)$ and assume in addition there exists a (piecewise smooth) timelike curve $\sigma:[0,1] \rightarrow M$ with $\sigma(0) =p$, $\sigma(1) =q$ which does not intersect $Conj_c(p)$. Then there exists a timelike geodesic from $p$ to $q$. 
\end{thm}
\noindent {\it Proof.} We shall construct a (piecewise smooth) {\it lift} of $\sigma$ through $\exp _p$ to $\mathcal{T}_p\cap \mathcal{D}_p$, that is, a piecewise smooth curve $\overline{\sigma}:[0,1] \rightarrow \mathcal{T}_p\cap \mathcal{D}_p$ with $\overline{\sigma}(0) = 0_p$ and $\exp_p \circ \overline{\sigma} = \sigma$. Since we then in particular have $\exp_p(\overline{\sigma}(1)) = q$ this establishes the conclusion. 

We apply here a slight modification of the main argument in the proof of \cite[Thm 2.11]{rein}. Since $0_p \notin \mathcal{S}_p$, by the Inverse Function theorem there exist open sets $U \ni 0_p$ in $\mathcal{D}_p$ and $V\ni p$ in $M$ such that $V=\exp_p(U)$ and $\exp_p|_{U}:U\rightarrow V$ is a diffeomorphism. By continuity of $\sigma$ there exists $0<\varepsilon <1$ with $\sigma[0,\varepsilon] \subset V$, so we set $\overline{\sigma}|_{[0,\varepsilon]}:= (\exp_p|_{U})^{-1} \circ \sigma |_{[0,\varepsilon]}$. 

The key observation here is that since $\exp_p \circ \overline{\sigma}|_{[0,\varepsilon]}$ is timelike by contruction, by \cite[Lemma 5.33]{oneill}, we have $\overline{\sigma}[0,\varepsilon]\subset \mathcal{T}_p$ (and indeed $\overline{\alpha}$ stays within a single timecone).
 
Let $0<\ell \leq 1$ be the lowest upper bound of all $t\in (0,1]$ for which $\overline{\sigma}|_{[0,t]}$ is well-defined, piecewise smooth and contained in $\mathcal{T}_p\cap \mathcal{D}_p$. 

Since $\exp_p\circ \overline{\sigma}|_{[0,\ell)} = \sigma |_{[0,\ell)}$ and $\overline{\sigma}(0)=0_p$, the CCP implies 
there exists some sequence $(t_k) \subset [0,\ell)$ with $t_k\rightarrow \ell$ for which $\overline{\sigma}(t_k) \rightarrow \hat{x}$ for some $\hat{x}\in C_p =\mathcal{C}_p \cap \mathcal{D}_p$. Put $\overline{\sigma}(\ell):= \hat{x}$. Since $\exp_p(\hat{x}) = \sigma(\ell)\notin Conj_c(p)$, $\hat{x} \notin \mathcal{S}_p$, again by continuity and the Inverse Mapping theorem we conclude that $\overline{\sigma}_{[0,\ell]}$ is well-defined, and again contained in $\mathcal{T}_p\setminus \mathcal{S}_p$ by \cite[Lemma 5.33]{oneill}. In addition, if $\ell <1$ we could again extend $\overline{\sigma}$ via a local inverse for $\exp_p$ contradicting the definition of $\ell$. 
\qcd

{\it The hypothesis of causal continuation cannot be removed in Theorem \ref{thmclave3}}. To see this, just consider the flat Lorentzian manifold $(M:=\mathbb{R}^2\setminus \{(1,0)\},-dt^2+dx^2)$, $p=(0,0), q=(2,0)$. Then $q\in I(p)$, but there is no timelike geodesic connecting them. Indeed, 
$$\mathcal{D}_p \simeq \mathbb{R}^2 \setminus \{(t,0)\, : \, t\geq 1\}.$$
Given any timelike curve $\sigma:[0,1] \rightarrow M$ from $p$ to $q$, its portion $\sigma|_{[0,1)}$ admits a lift to $C_p$ through $\exp _p$, but cannot be extended in $\mathcal{D}_p$.

\begin{cor}\label{timelikeloop}
Let $(M,g)$ be a Lorentzian manifold. If for some $p\in M$, $Conj_c(p)=\emptyset$ and $\exp_p$ has the CCP, then for any $q \in I(p)$ there exists a timelike geodesic from $p$ to $q$. In particular, if $p=q$ then there exists a timelike loop at $p$. 
\end{cor}
\qcd

\begin{rem}\label{remflaherty}
Suppose that for any $p \in M$ and for any {\it timelike plane} $\Pi_p \subset T_pM$, the sectional curvatures $K(\Pi_p)\geq 0$. Then for any $p \in M$, we have $Conj_c(p) =\emptyset$ by \cite[Prop. 2.1]{Flaherty1} (cf. also \cite[Prop. 11.13]{beem}). 
\end{rem}

Since Lorentzian coverings are in particular local isometries, Corollary \ref{timelikeloop}, Remark \ref{remflaherty}, and Proposition \ref{prop2} immediately yield 
\begin{cor}\label{coolcor} Let $(M,g)$ be a Lorenztian manifold with non-negative sectional curvatures on timelike planes. If $(M,g)$ admits a Lorentzian covering $\phi:(\hat{M},\hat{g})\rightarrow (M,g)$ such that $(\hat{M},\hat{g})$ is causally pseudoconvex and disprisoning, then for any $p,q \in M$ with $q\in J(p)$ there exists a causal geodesic connecting $p$ and $q$. 
\end{cor}
\qcd 

\section{A novel Lorentzian Hadamard-Cartan theorem}\label{sec4}

If a Lorentzian manifold $(M,g)$ is causally pseudoconvex and disprisoning (or has a Lorentzian covering that is) - which as seen in Remark \ref{ccpweaker} is strictly stronger than just having the CCP for the exponential maps - then much more can be said about causal geodesic connectedness. The results in this section are more interesting when applied to {\it spacetimes}, i.e., connected time-oriented Lorentz manifolds; accordingly, henceforth we assume that $(M,g)$ is time-oriented.   

Our first result here is the following extended version of the so-called {\it Lorentzian Hadamard-Cartan Theorem} (cf. \cite[Thm. 11.16]{beem}) which does not assume either global hyperbolicity (cf. Remark \ref{GH1}) or future 1-connectedness. 

\begin{thm}[Lorentzian Hadamard-Cartan]\label{lorentzianhadamard}
Let $(M,g)$ be a spacetime. Assume that $(M,g)$ is causally pseudoconvex and disprisoning. Let $p \in M$ be such that $Conj_c(p) =\emptyset$. Then the following hold. 
\begin{itemize}
    \item[a)] $\exp_p|_{\mathcal{T}^+_p \cap \mathcal{D}_p}: \mathcal{T}^+_p \cap \mathcal{D}_p \rightarrow I^+(p) $ is a diffeomorphism; in particular $I^+(p)$ is diffeomorphic to $\mathbb{R}^m$ and for any $q\in I^+(p)$ there exists a unique (up to reparametrization) future-directed timelike geodesic from $p$ to $q$. 
    \item[b)] If $q\in J^+(p)\setminus \{p\}$, there exist at most finitely many future-directed null geodesic segments (up to reparametrization) from $p$ to $q$.
    \item[c)] $J^+(p)$ is closed. 
\end{itemize}
\end{thm}
\noindent {\it Proof.} First, write $C_p= \mathcal{C}_p \cap \mathcal{D}_p$, $\tau _p:= \mathcal{T}_p \cap \mathcal{D}_p$, and $\tau^{+}_p:= \mathcal{T}^{+}_p\cap \mathcal{D}_p$, where $\mathcal{T}^{+}_p$ is the future timecone at $p$. 

By Corollary \ref{causalcorshort}, $\exp_p|_{C_p}$ is proper, and since $\mathcal{S}_p \cap \mathcal{C}_p = \emptyset$ by the hypothesis that $Conj_c(p) =\emptyset$, we also have that $\exp_p$ is a local diffeomorphism around each $v \in \mathcal{C}_p$. 

Let $q\in M$. Since $(\exp_p|_{C_p})^{-1}(q)$ is compact in $C_p$, if it were infinite one would have an accumulation point $v\in C_p$, which would contradict the local injectivity of $\exp_p$ around this point. Therefore, $(\exp_p|_{C_p})^{-1}(q)$ is finite (maybe empty). This already establishes $(b)$. 

By Prop. \ref{prop2} $\exp_p$ has the CCP, and hence by (the proof of) Theorem \ref{thmclave3} $\exp_p|_{\tau_p}: \tau_p \rightarrow I(p)$ is an onto map. Indeed, this proof can be easily adapted to show that $\exp_p|_{\tau^+_p}: \tau^+_p \rightarrow I^+(p)$ is an onto map. Since its restriction to the open set $\tau^+_p$ is still a proper map, and since it is a local diffeomorphism, by Prop. \ref{resultinbowder} it is a covering map (with finite fibers). 

However, $\tau^+_p$ is homeomorphic to $\mathbb{R}^m$, so $\exp_p|_{\tau^+_p}$ is a universal covering. By \cite[Lemma 8]{beemparker} $\pi_1(I^+(p))$ is either trivial or infinite, but we saw that the cardinality of the fibers of $\exp_p|_{\tau^+_p}$ - which equals that of $\pi_1(I^+(p))$ for an universal cover - is finite. We conclude that $\exp_p|_{\tau^+_p}: \tau^+_p \rightarrow I^+(p)$ is a trivial covering, and therefore a diffeomorphsim. This establishes $(a)$. 

Finally, to prove $(c)$ we recall that 
$$\overline{J^+(p)} = \overline{I^+(p)}.$$
Given $q\in \overline{J^{+}(p)}$, therefore, we can pick a sequence $(q_k)\subset I^+(p)$ converging in $M$ ro $q$. By item $(a)$, for each $k \in \mathbb{N}$, there exists a unique $v_k \in \tau^{+}_p\subset C_p$ such that $\exp_p(v_k)=q_k$. By the properness of $\exp_{p}|_{C_p}$ we conclude that up to passing to a subsequence we can assume that $v_k\rightarrow v_0$ for some $v_0 \in \overline{\tau^+_p}$, whence we conclude that $q\equiv \exp_p(v_0)$; it follows that $q\in J^+(p)$, which completes the proof. 
\qcd
\begin{rem}\label{semifinal}
Clearly, by time duality Theorem \ref{lorentzianhadamard} can also be analogously stated (and remains valid) for {\it pasts} of points.
\end{rem}
Recall that a spacetime $(M,g)$ is said to be {\it causally simple} if it is {\it causal}, i.e., has no closed causal curves, and for any $p\in M$ $J^{\pm}(p)$ are closed \cite{minguzzi}. 

We end with the following variant of both \cite{krolak} and \cite{Flaherty2}.
\begin{cor}\label{causalladderplace}
Let $(M,g)$ be a causal, causally pseudoconvex and disprisoning spacetime. Assume that for any $p \in M$, we have $Conj_c(p) =\emptyset$. Then $(M,g)$ is causally simple. If in addition $(M,g)$ is causally geodesically complete, then it is globally hyperbolic. 
\end{cor}
\noindent {\it Proof.}  By taking Theorem \ref{lorentzianhadamard}(c) and Remark \ref{semifinal} into account the simple causality immediately follows. 

 For the other claim, by taking the main theorem of \cite{Flaherty2} into account, what remains to be shown is that $(M,g)$ is {\it timelike $1$-connected}, i.e., given any $p,q\in  M$ with $q\in I^{+}(p)$, any two future-directed timelike curve segments from $p$ to $q$ are fixed-endpoint-homotopic through future-directed timelike curve segments from $p$ to $q$. 
 
  By Theorem \ref{lorentzianhadamard}, it suffices to show that any future-directed timelike curve $\alpha:[0,1]\rightarrow M$ with $\alpha(0)=p,\alpha(1)=q$ is fixed-endpoint-homotopic through timelike curves to the unique future-directed timelike geodesic segment $\eta:[0,1]\rightarrow M$ from $p$ to $q$. 
  
 Theorem \ref{lorentzianhadamard}(a) implies that there exists a unique piecewise smooth curve segment $\beta:[0,1]\rightarrow T_pM$ with $\beta(0)=0_p$ and $\beta(0,1]\subset  \mathcal{T}_p^{+}\cap \mathcal{D}$ such that $\exp_p\circ \beta = \alpha$. We can also write 
 $$\sigma(t)=\exp_p(t\cdot \beta(1)), \quad t\in [0,1].$$
Thus, define the manifestly continuous map $\overline{H}:[0,1] 
^2 \rightarrow (\mathcal{T}_p^{+}\cap \mathcal{D})\cup \{0_p\} $ by 
\begin{equation}
 \overline{H}(s,t):=   \left\{\begin{array}{cc}
        (t/s)\cdot \beta(s), &  \hbox{if } 0\leq t \leq s,\\
       \beta(t), & \hbox{if }s \leq t\leq 1,
    \end{array}\right.,
\end{equation}
when $0<s\leq 1$; and $\overline{H}(0,t):= \beta(t)$ for all $t\in [0,1]$. Defining $H:= \exp _p\circ \overline{H}$, one easily checks this is the desired fixed-endpoint homotopy through timelike curves from $\alpha$ to $\sigma$, so the proof is complete. 


 

\qcd


\end{document}